\documentclass[a4paper,11pt,reqno]{amsart}
\usepackage[utf8]{inputenc}
\usepackage{mathtools}
\usepackage{amssymb}
\usepackage{amsfonts}
\usepackage{amsthm}
\usepackage{bbm}
\usepackage{fullpage}
\usepackage{color}
\usepackage{hyperref}
\usepackage{stmaryrd}
\usepackage{enumitem}
\usepackage{accents}
\mathtoolsset{showonlyrefs}
\newcommand{\N}{\mathbb N}

\newcommand{\E}{\mathbb E}

\newcommand{\ind}{\mathbbm 1}

\newtheorem{theorem}{Theorem}
\newtheorem{lemma}[theorem]{Lemma}
\newtheorem{proposition}[theorem]{Proposition}
\newtheorem{corollary}[theorem]{Corollary}
\newtheorem*{conjecture}{Conjecture}

\newtheorem{remark}[theorem]{Remark}

\theoremstyle{definition}

\newcommand{\dd}{\text{d}}

\newcommand{\p}{{p^*}}

\newcounter{constants}

\newcommand{\esssup}{\operatorname{ess\,sup}}

\newcommand{\bbE}{{\ensuremath{\mathbb E}} }

\newcommand{\bbP}{{\ensuremath{\mathbb P}} }

\newcommand{\bbR}{{\ensuremath{\mathbb R}} }

\newcommand{\bbZ}{{\ensuremath{\mathbb Z}} }

\renewcommand{\epsilon}{\varepsilon}

\newcommand{\gep}{\varepsilon}       % \ge already exists...

\newcommand{\go}{\omega}

\newcommand{\gO}{\Omega}
\newcommand{\gl}{\lambda}

\newcommand{\cF}{{\ensuremath{\mathcal F}} }

\newcommand{\cC}{{\ensuremath{\mathcal C}} }

\newcommand{\cZ}{{\ensuremath{\mathcal Z}} }

\newcommand{\bP}{{\ensuremath{\mathbf P}} }

\newcommand{\bE}{{\ensuremath{\mathbf E}} }

\newcommand{\lint}{\llbracket}
\newcommand{\rint}{\rrbracket}

\newcommand{\uniontwo}[2]{\bigcup_{\substack{#1 \\ #2}}} % union with 2 lines
\definecolor{darkred}{rgb}{0.7,0.1,0.1}
\renewcommand{\hat}{\widehat}

\title{Local limit Theorem for Directed polymers in a  random environment in the full weak disorder regime}

\author{Stefan Junk}
\address{Gakushuin University, 1-5-1 Mejiro, Toshima-ku, Tokyo 171-8588 Japan}
\email{sjunk@math.gakushuin.ac.jp}
\author{Hubert Lacoin}

\address{IMPA, Estrada Dona Castorina 110,
Rio de Janeiro RJ-22460-320- Brasil \\
Università Roma Tre,
Largo San Murialdo 1, 00146 Roma - Italia}
\email{lacoin@impa.br}

\begin{document}

\maketitle

\begin{abstract}
For directed polymers in random environment, we provide a short and mostly self-contained proof of the local limit theorem in the full weak disorder phase, including the critical inverse temperature $\beta_c$. That is, the point-to-point partition function is well approximated by the product of two independent point-to-line partition functions (one starting from each endpoint of the system) multiplied by the simple random walk transition kernel. The approximation is valid in the $L^p$ sense for an optimal range of $p$. The result is new at $\beta=\beta_c$ and for $p>2$, and in all cases the proof is significantly simpler, technically and conceptually, than what can be found in the literature. As a consequence of the result we establish an invariance principle for the polymer measure with fixed end-point, valid for all $\beta \in [0,\beta_c]$.
\end{abstract}

%\begin{abstract}
%We provide a short and mostly self contained proof that for every $\beta\in [0,\beta_c]$ -- where $\beta_c$ denote the critical inverse temperature -- the point-to-point partition function of the directed polymers in a random environment is well approximated by the product of two independent point-to-line partition function (one starting from each extremity of the system). The approximation is valid in the $L^p$ sense for an optimal range of $p$. The result is new when $\beta=\beta_c$, and for $\beta\in (0,\beta_c)$ the proof is significantly simpler (technically and conceptually) than
%what can be found in the literature. As a consequence of the result we establish an invariance principle for the polymer measure with constrained end-point, valid for all $\beta \in [0,\beta_c]$.
%\end{abstract}

\section{Introduction and main result}

Let $(\go_{k,x})_{n\in \bbZ,x\in \bbZ^d}$ be a collection of i.i.d.\ random variables which satisfy
\begin{equation}
 \forall \beta\in \bbR, \quad  \bbE\left[ e^{\beta\go_{1,0}-\gl(\beta)}\right]<\infty.
\end{equation}
For a fixed realization of $\go$ we define the polymer measure $P^{\beta,\go}_n$ as the following modification of the simple random walk measure $P$:
\begin{equation}
 \dd P^{\beta,\go}_{n}= \frac{1}{W^{\beta}_n}\exp\left( \sum_{k=1}^n \beta\go_{k,X_k}-\gl(\beta)\right) \dd P  \quad   \text{ where } W^{\beta}_n:=E\left[ \exp\left( \sum_{k=1}^n \beta\go_{k,X_k}-\gl(\beta)\right)\right].
\end{equation}
It is known  that for $d\ge 3$ there exists $\beta_c\in (0,\infty]$  such that
\begin{equation}
\forall \beta\in [0,\beta_c]\cap  [0,\infty) , \quad    \bbP\left( W^{\beta}_{\infty} >0 \right)=1.
\end{equation}
More precisely, the fact that  $\bbP( W^{\beta}_{\infty} >0)=1$ for $\beta$ small was established in \cite{B89,IS88}. In \cite{CY06} it was proved that 
$\bbP( W^{\beta}_{\infty} >0 )$ is monotone in $\beta$ (hence the existence of $\beta_c$) and finally in \cite{JL_ASENS,JL_ORB} it was shown that $\bbP( W^{\beta_c}_{\infty} >0)=1$ in the case where $\beta_c$ is finite (when $d=1,2$ we have $\beta_c=0$, see \cite{CH02,CSY03}).
Note that it is known that $\beta_c<\infty$ if $\bbP(\go_{1,0}=\operatorname{ess\, sup}\go_{1,0})<\frac 1{2d}$. Since we only consider finite values of $\beta$, in the remainder of the paper we abuse notation and write $[0,\beta_c]$ for $[0,\beta_c]\cap [0,\infty)$.
%Following a terminology introduced in  \cite{J22}, we set for  $\beta\in (0,\beta_c]$  
%\begin{equation}
% p^*(\beta):= \sup\left\{ p>1 \ : \  \bbE\left[ (W^{\beta}_n)^p \right]<\infty \right\}.
%\end{equation}

\medskip
It has been established in \cite{CY06}, that for $\beta\in [0,\beta_c]$, the scaling limit of $(X_k)_{k=1}^n$ under $P^{\beta,\go}_{n}$ is the same as the for the simple random walk.
In particular for any continous bounded function on $\bbR^d$ we have 
\begin{equation}\label{weakclt}
 \lim_{n\to \infty }E^{\beta,\go}_{n}\left[\varphi\left(\frac{X_n}{\sqrt{n}}\right) \right]= \left(\frac{d}{2\pi}\right)^{d/2} \int_{\bbR^d} \varphi(u) e^{-\frac{du^2}{2}} \dd u.
\end{equation}
We refer to \cite{AZ96,B89}, for earlier versions of  this invariance principle and central limit Theorem proved under more restrictive assumption.
A natural follow-up question is whether one can derive a local version of this central limit theorem.

\medskip

A local CLT should take the form of a \textit{factorization} of the point to point partition function, which we define by 
\begin{equation}\label{ptop}
 \hat W^{\beta}_n(x):= E\left[ \exp\left( \sum_{k=1}^n \beta\go_{k,X_k}-\gl(\beta)\right)\ind_{\{X_n=x\}}\right].
\end{equation}
Given a space time point $(n,x)$  and $m\ge 0$ we define
\begin{equation}\label{backarow}
 \overleftarrow{W}^{\beta}_{m}(n,x):=E\left[ \exp\left( \sum_{k=0}^m \beta\go_{n-k,x-X_k}-\gl(\beta)\right)\right],
\end{equation}
This is the partition function of a polymer model of length $m$ starting from $(n,x)$ and running backward in time. Note that, for reason that will become obvious, the environment at the starting point is included in $\overleftarrow{W}^{\beta}_{m}(n,x)$. We define $\overleftarrow{W}^{\beta}_{\infty}(n,x)=\lim_{m\to \infty}\overleftarrow{W}^{\beta}_{m}(n,x)$.

\medskip
From a local CLT, we expect that when $x$ is of order $\sqrt{n}$ and of the same parity as $n$, then as $n\to \infty$ one should have
\begin{equation}\label{approxx}\begin{split}
 \hat W^{\beta}_n(x)&= P(X_{n}=x) \left(W^{\beta}_{\infty} \times  \overleftarrow{W}^{\beta}_{\infty}(n,x) +o(1)\right)\\
 &= 2\left(\frac{d}{2\pi n}\right)^{d/2}  e^{-\frac{d|x|^2}{2n}}\left(W^{\beta}_{\infty} \times  \overleftarrow{W}^{\beta}_{\infty}(n,x) +o(1)\right)
\end{split}\end{equation}
(the second line is merely the local CLT of the unperturbed measure $P$).
The influence of the environment is felt only at the two tips of the system (near the starting point and the endpoint) while the environment in the bulk homogenizes, so asymptotically the partition function factorizes. In the limit $n\to\infty$, the point-to-point partition function is well approximated by the product of two infinite volume point-to-line partition functions. Moreover, for the polymer measure $P_n^{\beta,\go}$  the denominator cancels and we expect
\begin{align*}
	P_n^{\beta,\go}(X_n=x)=\frac{\hat W_n^\beta(x)}{W_n^\beta}=P(X_n=x)\overleftarrow W^\beta_\infty(n,x)(1+o(1)).
\end{align*}

\medskip
\subsection{The main result}

Our main achievement in this paper it a relatively short and simple proof   that
\eqref{approxx} holds in the $L^p$ sense, for the largest possible range of $\beta$ and $p$. To make this precise, following a terminology introduced in  \cite{J22}, we set for  $\beta\in (0,\beta_c]$
\begin{equation}
 p^*(\beta):= \sup\left\{ p>1 \colon \sup_{n\ge 1}  \bbE\left[ (W^{\beta}_n)^p \right]<\infty \right\}.
\end{equation}

\begin{theorem}\label{factoth}
For any $\beta\in[0, \beta_c]$, $p\in [1,p^*(\beta))$ and sequences $(\ell_n)_n$ and $(m_n)_n$ such that $(\ell_n+ m_n+1)\le n$ and $\lim_{n\to \infty} \ell_n \wedge m_n=\infty$ we have

\begin{equation}\label{facto}
\lim_{n\to \infty}n^{\frac{pd}{2}} \sup_{x\in \bbZ^d}\bbE\left[ \left|\hat W^{\beta}_{n}(x)- W^{\beta}_{\ell_n}\times \overleftarrow{W}^{\beta}_{m_n}(n,x)\times P(X_{n}=x)\right|^p\right]=0.
\end{equation}
\end{theorem}

\begin{remark}
In the above result, the restriction $(\ell_n+m_n+1)\le n$ is taken to ensure that $W^{\beta}_{\ell_n}$ and  $\overleftarrow{W}^{\beta}_{m_n}(n,x)$ are independent. Without this restriction we have no guarantee that the product $W^{\beta}_{\ell_n}\times \overleftarrow{W}^{\beta}_{m_n}(n,x)$ remains bounded in $L^p$.
\end{remark}

\subsection{Comparison with the existing literature}

Earlier versions of Theorem~\ref{factoth} have been proved with more restrictive conditions. Let us review here earlier results and underline where the novelty of the present paper lies.
The factorization formula  \eqref{facto} in Theorem \ref{factoth} was proved first in \cite{sinai95} (see also \cite{V06} for a more recent revisit)  for $p=2$ under the assumption that $\beta< \beta_2$, where
\begin{equation}\label{l2thresh}
\beta_2:= \sup\left\{ \beta>0 : \sup_n\bbE\left[ (W^{\beta}_n)^2\right]<\infty\right\}.
\end{equation}
The proofs in \cite{sinai95,V06} take advantage of the explicit expression for the second moment of the partition function: it can be written in terms of the Laplace transform of the number of collisions between two independent random walks and allows explicit calculations.

\medskip

However, in \cite{BT10,BS10,BS11} it was proved that the $L^2$-integrability condition imposed by $\beta_2$ is strictly stronger than weak disorder, i.e.,  $\beta_2<\beta_c$. In particular, the results obtained in \cite{sinai95,V06} leave open the question of whether some form of factorization holds when $\beta$ is in the interval $[\beta_2,\beta_c]$.

\medskip

The first result establishing \eqref{facto} beyond the $L^2$-regime was obtained recently in \cite[Theorem~1.1]{J_JEMS}, where it is proved with two additional restrictions:
\begin{itemize}
 \item[(a)] $p^*(\beta)>1+\frac{2}{d}$
 \item [(b)] $p<2\wedge p^*(\beta)$.
\end{itemize}
It was also shown \cite[Theorem~1.2]{J_JEMS} that (a) is equivalent to $\beta\in[0,\beta_c)$  under a restrictivte assumption on the environment, which will be removed in an upcoming work.
The reason why restriction (a) is present is because some key computations in the proof require
\begin{equation}
 \sum_{n\ge 1}\sum_{x\in \bbZ^d} P(X_n=x)^{p}<\infty.
\end{equation}
On the other hand, the restriction to $p< 2$ plays a less crucial role and we believe that it can possibly be waived at the cost of additional technical complications
(in the same spirit as what is done in Section \ref{p_larger_2} of the present paper).

The novelty of the present work is thus to give a proof of Theorem \ref{facto} which is simpler than the one presented in \cite{J_JEMS} and extends to the case $\beta=\beta_c$ and $p\ge 2$.
Our approach borrows and adapts many of the important ideas developed in \cite{J_JEMS},
but adds a key new input that allows to bypass many difficult computation: a reccursive approach to prove $L^p$-tightness of
the normalized point-to-point partition function $\hat W_n(x)/P(X_n=x)$.

\subsection{Discussing the case of endpoints  on a superdiffusive scale}

When $x\gg \sqrt{n}$, the result as written in \eqref{facto} only tells us that
$\bbE[(\hat W_n(x))^p]=o(n^{-\frac{dp}{2}})$, but does not provide a factorization formula in the spirit of \eqref{approxx}. Let us discuss below what is expected beyond the diffusive scale.

\subsubsection*{An easy but useful bound}
With very little effort we can obtain a more quantitative bound.

\begin{corollary}
For every $\beta\leq \beta_c$ and $p<p^*(\beta)$, there exists $C(\beta,p)>0$ such that for every $x\in \bbZ^d$
\begin{equation}\label{desired}
	\bbE\left[(\hat W^{\beta}_{n}(x))^p\right] \le C(\beta,p) n^{-\frac{dp}{2}} e^{-\frac{|x|^2}{C(\beta,p)n}}.
\end{equation}
Moreover, $C(\beta,p)$ can be chosen non-decreasing in $\beta$ and $p$.
\end{corollary}
\begin{proof}
We observe that \eqref{facto} implies that for any $p<p^*(\beta)$ we have
\begin{equation}\label{boundedd}
 \limsup_{n\to \infty}n^{\frac{pd}{2}} \sup_{x\in \bbZ^d}\bbE\left[ \left|\hat W^{\beta}_{n}(x)\right|^p\right]= \bbE\left[ (W^{\beta}_\infty)^p \right] \bbE\left[ (\overleftarrow{W}^{\beta}_\infty)^p \right] \limsup_{n\to \infty}n^{\frac{pd}{2}} \sup_{x\in \bbZ^d} P(X_{n}=x)^p<\infty.
\end{equation}
Furthermore, the following quantity is monotone in $\beta$ and $p$
\begin{equation}
S(\beta,p):= \sup_{n\ge 1} \sup_{x\in \bbZ^d} n^{\frac{d}{2}} \bbE\left[ \left|\hat W^{\beta}_{n}(x)\right|^p\right]^\frac{1}{p}.
\end{equation}
By \eqref{boundedd}, $S(\beta,p)$ is finite for $1\le p<p^*(\beta)$.
Monotonicity in $p$ follows from H\"older and that in $\beta$ from convex comparison (see e.g. \cite[Lemma 3.3]{CY06} or \cite[Lemma 5.2]{L25bb}).
Given $p<p^*(\beta)$ and $p'\in(p,p^*(\beta))$, we can use H\"older's inequality to obtain
\begin{equation}\label{trey1}\begin{split}
 n^{\frac{dp}{2}} \bbE\left[(\hat W^{\beta}_{n}(x))^p\right]&\le
 n^{\frac{dp}{2}}\bbE\left[(\hat W^{\beta}_{n}(x))^{p'}\right]^{\frac{p-1}{p'-1}}
  P(X_n=x)^{\frac{p'-p}{p'-1}}\\
  &\le S(\beta,p')^\frac{p'(p-1)}{p'-1} \left( n^{\frac{d}{2}} P(X_n=x) \right)^{\frac{p'-p}{p'-1}}.
\end{split}\end{equation}
The usual upper bound for the simple random walk yields $n^{\frac{d}{2}} P(X_n=x)\le C e^{-\frac{|x|^2}{4n}}$ for some $C=C(d)\ge 1$ and hence we obtain (with the same $C$)
\begin{equation}
n^{\frac{dp}{2}} \bbE\left[(\hat W^{\beta}_{n}(x))^p\right]
\le \left(S(\beta,p')\vee C\right)^p e^{-\frac{p'-p}{p'-1} \frac{|x|^2}{4n}}.
\end{equation}

We therefore set
\begin{align*}
C(\beta,p):=\inf_{p<p'<p^*(\beta)}\left\{\left(S(\beta,p')\vee C\right)^p\vee\frac{4(p'-1)}{p'-p}\right\}<\infty.
\end{align*}
The inequality \eqref{desired} holds for every choice of $p<p'<\p(\beta)$, hence it remains valid upon taking the infimum. It remains to show that $C(\beta,p)$ is increasing in both arguments: note that for fixed $p'$, both terms are increasing in both arguments. Moreover, since $\p$ is decreasing, the admissible interval for $p'$ shrinks when either $\beta$ or $p$ increases, which yields the desired monotonicity.
\end{proof}

\subsubsection*{Factorization beyond the diffusive scale}
While Theorem~\ref{factoth} only establishes factorization of the point-to-point partition function (in the sense of \eqref{approxx}) for $|x|=O(\sqrt{n})$, the factorization \eqref{approxx} should hold as soon as $|x|=o(n)$. While we present this fact as a conjecture, let us mention that it seems very plausible that the techniques developed in \cite{J_JEMS} are sufficient to prove it.

\begin{conjecture}
 Let $\beta\in (0,\beta_c)$, $p<p^*(\beta)$, $x_n$, $\ell_n$ and $m_n$ sequences that satisfy
 $$P(X_n=x_n)>0, \quad \lim_{n\to \infty} |x_n|/n=0, \quad (\ell_n+ m_n+1)\le n, \text{ and } \lim_{n\to \infty} \ell_n \wedge m_n=\infty.$$
We have
\begin{equation}\label{superfacto}
\lim_{n\to \infty}\bbE\left[ \left| \frac{\hat W^{\beta}_{n}(x_n)}{P(X_n=x_n)}- W^{\beta}_{\ell_n}\times \overleftarrow{W}^{\beta}_{m_n}(n,x_n)\right|^p\right]=0.
\end{equation}
\end{conjecture}
Theorem~\ref{factoth} proved in this paper only implies that \eqref{superfacto} holds under the stronger assumption
$\limsup_{n\to \infty} |x_n| n^{-1/2}<\infty$.
Note that the condition $|x_n|=o(n)$ is necessary. If $\lim_{n\to \infty} x_n= \alpha \in (-1,1)^d$,
then  conditionning on the event $\{X_n=x_n\}$ induces a tilt on the walk $X$. The distribution of $P(X\in \cdot \ | \ X_n=x_n)$ converges to the law $Q_{\alpha}$ of a random walk with i.i.d.\ increments distributed as 
\begin{equation}
 Q_{\alpha}(X_1 \in A)= \frac{E\left[ e^{\langle\theta, X_1\rangle} \ind_A\right]}{E\left[ e^{\langle \theta,X_1\rangle} \right]} \quad \text{ where $\theta\in \bbR^d$ satisfies} \quad    \alpha= \frac{E\left[ X_1 e^{\langle\theta,X_1\rangle} \right]}{E\left[ e^{\langle\theta,X_1\rangle} \right]}.
\end{equation}

In that case,  the analogue of \eqref{superfacto} should still holds but with $W^{\beta}_{\ell_n}$ and $\overleftarrow{W}^{\beta}_{m_n}(n,x_n)$  replaced by the  partition functions obtained with an average over $Q_{\alpha}$, and  $p^*(\beta)$ has to be replaced with $p^*(\beta,\alpha)$, defined just as $p^*(\beta)$ but with $P$ replaced by $Q_{\alpha}$
(see e.g. \cite[Section~1.6.3]{J_JEMS}). Note also that it has been proved \cite[Corollary~1.9]{J_JEMS} that $\liminf_{\alpha\to 0}\p(\beta,\alpha)=\p(\beta)$ for $\beta<\beta_c$, so that the assumption $p<p^*(\beta)$ implies that $p<p^*(\beta,\alpha)$  for  sufficiently small $\alpha$. Hence when $\beta<\beta_c$, we may expect that -- with appropriate modifications -- the factorization \eqref{superfacto} extends to the range  $|x_n|\le \delta(p) n$.

\medskip

When $\beta=\beta_c$, while we believe that while \eqref{superfacto} should remain true well beyond the diffusive window, it should require a stronger condition than $|x_n|=o(n)$ to hold true. This is because even when the tilt induced by the conditionning on $\{ X_n=x_n\}$ is vanishing, the induced decrease of entropy might still be sufficient to push the system into supercriticality, in analogy with what happens when considering sequences $\beta_n\downarrow\beta_c$ approaching criticality from above (see \cite[Proposition 2.10]{JL_in_prep}). In particular we believe that the critical point $\beta_c(\alpha)$ for the directed polymer constructed with $P$ replaced by $Q^{\alpha}$ satisfy $\beta_c(\alpha)<\beta_c$ for every $\alpha\ne 0$.
\medskip

Finally we note that \eqref{facto} is proved for $p=2$ and $\beta<\beta_2$ in \cite[Theorem~2]{HKL23} under the stricter condition $|x_n|\le n^{\sigma}$ for some $\sigma<1$. When $\beta<\beta_c$, it is established in \cite[Corollary 1.6]{J_JEMS} that \eqref{superfacto} is valid for $p\in(1+\frac 2d,2\wedge\p(\beta))$ up to $|x_n|\le n^{\sigma}$ for $\sigma<3/4$.
 We strongly believe that the techniques developed in these references can be extended to treat any case with $|x_n|=o(n)$.

\subsection{A conditionned invariance principle}

Let us now present a consequence of Theorem~\ref{factoth} on the distribution of the polymer trajectory under the \textit{conditioned} polymer measure.
We extend the random walk $(X_k)_{k\in \bbZ_+}$ to noninteger times by affine extension by setting, for $s\in [0,\infty]$, 
$$X_s= (s+1 -\lceil s\rceil) X_{\lceil s\rceil} +  (\lceil s \rceil -s) X_{\lceil s\rceil-1}$$ 
and define the diffusively rescaled trajectory
$(X^{(n)}_t)_{t\in [0,1]}=\left( X_{ n t}/\sqrt{n}\right)_{t\in [0,1]}$.

\medskip

Note that $X^{(n)}$ is a random process taking value in $\cC:=C([0,1],\bbR^d)$, the set of continuous functions from $[0,1]$ to $\bbR^d$. We consider $\cC$ equiped with the topology of uniform convergence. It has been established in \cite{CY06} that $P^{\beta,\go}_{n}\left( X^{(n)}\in \cdot \right)$ converges in probability to the Wiener measure in weak disorder.
This means that, for any bounded continuous function $\varphi\colon \cC\to\bbR$ we have
\begin{equation}\label{invprincip}
\lim_{n\to \infty} E^{\beta,\go}_{n}\left[ \varphi( X^{(n)}) \right] = \bE\left[ \varphi( (B_t)_{t\in [0,1]}) \right],
\end{equation}
where $B$ is a Brownian Motion with diffusivity matrix $\frac{1}{d}I_d$ ($I_d$ is the identity matrix in $\bbR^d$).
An alternative shorter proof of \eqref{invprincip} is presented in \cite{L25}, which relies mostly on exploiting the convergence $\lim_{n\to \infty}W^{\beta}_n= W^{\beta}_\infty$ in $L^1$. Using the same ideas and the factorization \eqref{facto}, we can prove an analogue of \eqref{invprincip} for the polymer trajectory with a fixed endpoint.

\medskip

 Let us recall that the simple random walk with fixed endpoint converges in distribution to a \textit{Brownian Bridge}, that is to say, if $\lim_{n\to \infty} x_n/\sqrt{n}= {\bf x}\in \bbR^d$, then for any continous  function $\varphi\colon \cC\to [0,1]$

\begin{equation}\label{condinvprincip}
	\lim_{n\to \infty} E\left[ \varphi\left((X^{(n)}_t)_{t\in[0,1]}\right) \ \middle| \ X_n=x_n \right]= \bE\left[ \varphi\left( (B^{({\bf x})}_t)_{t\in [0,1]} \right) \right]
\end{equation}
where $\bP$ is the distribution of $(B_t)_{t\in [0,1]}$ as in \eqref{invprincip} and $B^{({\bf x})}$, defined by
\begin{align}\label{bb}B^{({\bf x})}_t:=B_t- tB_1+ t{\bf x}\end{align}
is a Brownian bridge with diffusivity matrix $I_d/d$ and endpoint ${\bf x}$. The convergence of finite dimensional marginals can be obtained using a local CLT (see for instance \cite[Chapter 2]{LL10}) and tightness can be obtained via moment estimates for the increments.

\begin{theorem}\label{condinvpr}
	Given  $\beta\in(0, \beta_c]$, a sequence $(x_n)_{n\in\N}$ such that $P(X_n=x_n)>0$ and $\lim_{n\to \infty} x_n/\sqrt{n}= {\bf x}$, a continous function $\varphi\colon \cC \to [0,1]$ and  $B^{(\bf x)}$ defined as in \eqref{bb}, the following   holds
\begin{equation}
\lim_{n\to \infty} E^{\beta,\go}_n\left[ \varphi(X^{(n)}_t) \ \middle| \ X_n=x_n \right]= \bE\left[ \varphi( (B^{(\bf x)}_t)_{t\in [0,1]} )  \right].
\end{equation}

\end{theorem}

\subsection{Notation and paper organization}

We let $W^{\beta}_n(\cdot)$ denote the unnormalized polymer measure and, given $x\in \bbZ^d$, $\hat W^{\beta}_n(\cdot,x)$ its point-to-point analogue by setting for $g:\gO\to \bbR$ measurable
\begin{equation}\begin{split}\label{Wasmeasure}
  W^{\beta}_n(g)&:=\bbE\left[ e^{\sum_{k=1}^n \beta \go_{k,X_k}-\gl(\beta)}g(X) \right],\\
    \hat W^{\beta}_n(g,x)&:=\bbE\left[ e^{\sum_{k=1}^n \beta \go_{k,X_k}-\gl(\beta)}g(X)\ind_{\{ X_n=x\}} \right].
\end{split}\end{equation}
For $D\subset \gO$ measurable, we set $W^{\beta}_n(D):= W^{\beta}_n(\ind_D)$ and $\hat W^{\beta}_n(D,x):=\hat W^{\beta}_n(\ind_D,x)$.  

\medskip In the next section, we explain the outline of the main result and give the proof of Theorem~\ref{factoth} assuming two key estimates, Propositions~\ref{factoslice} and \ref{keystep}. These propositions are then proved in Sections~\ref{sec_factoslice} and \ref{sec_keystep}. Finally, Theorem~\ref{condinvpr} is proved in Section~\ref{sec_condinvpr}.

\section{The outline of  the proof of Theorem \ref{factoth}}

\subsection{The main steps}

The first step of the proof is perform the factorization for a partition function in which a macroscopic slice of environment is removed through partial annealing.
Given $1\le a<b\le n$ we define
\begin{equation}\label{defwaabb}
\hat W^{(a,b)}_{n}(x):= E\left[  e^{\sum_{k\in \lint1,a\rint\cup\lint b,n\rint }\beta\go_{k,X_k}-\gl(\beta)}\ind_{\{X_n=x\}}     \right].
\end{equation}
Letting $\cF_A:= \sigma\left( \go_{n,x}, n\in A, x\in \bbZ^d\right)$ for $A\subset \bbZ$ we have
\begin{equation}\label{annil}
 \hat W^{(a,b)}_{n}(x) = \bbE\left[  \hat W^{\beta}_{n}(x) \ \middle| \ \cF_{\lint 1, a\rint\cup \lint b,n\rint} \right].
\end{equation}
Decomposing the partition function according to the position of the walk at time $a$ and $b-1$ and using Markov's property we observe that $\hat W^{(a,b)}_{n}(x)$ can be expressed using a product of point-to-line partition function as follows (recall \eqref{backarow})
\begin{equation}\begin{split}
 \hat W^{(a,b)}_{n}(x)&= \sum_{y,z\in \bbZ^d} \hat W^{\beta}_a(y)P(X_{b-a-1}=z-y) \theta_{b-1,z}\hat W^{\beta}_{n-b+1}(x-z)\\
 &\le   W^{\beta}_{a}\times \overleftarrow{W}^{\beta}_{n-b}(n,x)    \sup_{w\in \bbZ^d} P(X_{b-a-1}=w)\\
 &\le C (b-a)^{-\frac{d}{2}} W^{\beta}_{a}\times \overleftarrow{W}^{\beta}_{n-b}(n,x).
\end{split}\end{equation}
The last line is a consequence of the local CLT for $X$.
In the second line we used the fact that
\begin{equation}
 \sum_{y\in \bbZ^d} \hat W^{\beta}_a(y)= W^{\beta}_{a}   \quad \text{ and } \quad
 \sum_{z\in \bbZ^d}\theta_{b-1,z}\hat W^{\beta}_{n-b+1}(x-z)= \overleftarrow{W}^{\beta}_{n-b}(x,n).
\end{equation}
For $p<p^*(\beta)$, this yields in particular
\begin{equation}\label{superrough}
 \bbE\left[  \hat W^{(a,b)}_{n}(x)^p \right]\le (b-a)^{-\frac{pd}{2}} \bbE\left[ (W^{\beta}_{a})^p \left(\overleftarrow{W}^{\beta}_{n-b}(n,x)\right)^p\right]\le C_p (b-a)^{-\frac{np}{2}}.
\end{equation}
When $b-a$ is of order $n$, this bound is of the correct order of magnitude.
Pushing the above computation further and using the CLT for weak disorder polymer displayed in  \eqref{weakclt} we can obtain a factorization result for $\hat W^{(a,b)}_{n}(x)$.
To state it, we fix the values of $a_n$ and $b_n$. Given ${\bf a}$ and ${\bf b}$ satisfying $1/4\le{\bf a}<{\bf b}\le 3/4$, let
\begin{equation}\label{defanbn}
  a_n:=\lceil {\bf a} n\rceil \quad  \text{ and }  \quad b_n:=\lfloor {\bf b} n\rfloor.
\end{equation}

\begin{proposition}\label{factoslice}
If $\beta\leq \beta_c$ and $p<p^*(\beta)$,
  then
 \begin{equation}\label{keyddd}
\lim_{n\to \infty} n^{dp/2}\sup_{x\in \bbZ^d}\bbE\left[ \left|\hat W^{(a_n,b_n)}_{n}(x) -  W^{\beta}_{a_n}\times \overleftarrow{W}^{\beta}_{n-b_n}(n,x) P(X_{n}=x) \right|^p\right]=0.
 \end{equation}
\end{proposition}
Note that in particular \eqref{keyddd} implies
\begin{equation}\label{lessrough}\begin{split}
\lim_{n\to \infty} n^{\frac{dp}{2}}\sup_{x\in \bbZ^d}\bbE\left[ \hat W^{(a_n,b_n)}_{n}(x)^p\right]&=
\lim_{n\to \infty} n^{\frac{dp}{2}}\sup_{x\in \bbZ^d}\bbE\left[ (W^{\beta}_{a_n})^p\times \overleftarrow{W}^{\beta}_{n-b_n}(n,x)^p P(X_{n}=x)^p\right]\\
&= 2^p\left(\frac d {2\pi}\right)^{\frac{pd}{2}} e^{\gl(\beta p)-p\gl(\beta)} \bbE\left[  (W^{\beta}_\infty)^p\right]^2.
\end{split}\end{equation}
In particular, for later applications it is important that the right-hand side of \eqref{lessrough} is independent of ${\bf a}$ and $\bf b$.
While the presentation here is quite different from the one found \cite{J_JEMS}, the idea of erasing the environment in a slice to permit factorization is one of the main ingredient in the proof of \cite[Theorem 1.1]{J_JEMS}. 
The main novelty in our approach lies in the second step: we use an inductive approach to show that $\hat W^{(a,b)}_{n}(x)$ is a good approximation of $\hat W_{n}(x)$. The core of the argument is the following result.

\begin{proposition}\label{keystep}
Recall \eqref{defanbn}.
Given $p<p^*(\beta)$, there exists an increasing function $\gep:[0,1]\to \bbR_+$ with $\lim_{\delta \to 0} \gep(\delta)=0$
such that for all $n$ sufficiently large
\begin{equation}\label{subtle}
\bbE[\hat W_{n}(x)^p]\le  \bbE\left[\hat W^{(a_n,b_n)}_{n}(x)^p \right]  + \gep \left( {\bf b-a} \right) \sup_{y\in \bbZ^d}\bbE\left[  \hat W_{n-b_n}(y)^p \right].
\end{equation}

\end{proposition}
Let us try to explain here the meaning of the inequality \eqref{subtle} by proving a simpler inequality, which is valid for all $n,m\ge 1$ and $p<p^*(\beta)$,
\begin{equation}\label{rough2}\bbE[\hat W_{n+m}(x)^p]\le  C_p \sup_{y\in \bbZ^d}\bbE\left[  \hat W_m(y)^p \right].
\end{equation}
This claim follows from the following, more generic statement on the partition functions restricted to an event $D$ (recall \eqref{Wasmeasure}).
\begin{lemma}\label{sympa}
 If $D\in \sigma( (X_k)_{k\in [0,n]})$  then
 \begin{equation}
   \bbE\left[ \hat W_{n+m}(D, x)^p\right] \le
\bbE\left[ W_n(D)^p\right] \sup_{y\in \bbZ^d} \bbE\left[ \hat W_m(y)^p\right]
 \end{equation}
\end{lemma}

\begin{proof}[Proof of Lemma \ref{sympa}]
We start by decomposing the partition function according to the value taken by $X_n$. Using the Markov property for the random walk at time $n$ we obtain
\begin{equation}
  \hat W_{n+m}(D, x)=\sum_{y\in \bbZ^d}  \hat W_n(D,y) \theta_{n,y} \hat W_{m}(x-y)
\end{equation}
Then using Jensen's inequality for the probability $\hat W_n(D,\cdot)/W_n(D)$ (on $\bbZ^d$) we obtain
\begin{equation}
\left(\sum_{y\in \bbZ^d} \hat W_n(D,y) \theta_{n,y} \hat W_m(x-y)\right)^p
\le W_n(D)^{p-1} \sum_{y\in \bbZ^d}  \hat W_m(D,y)\left(\theta_{n,y} \hat W_{m}(x-y)\right)^p.
\end{equation}
From this we obtain (using independence and translation invariance) that
\begin{equation}\begin{split}
 \bbE\left[  \hat W_{n+m}(D,x)^p \ \middle| \ \cF_n  \right]&=  W_n(D)^{p-1} \sum_{y\in \bbZ^d} \hat W_n(D,y) \bbE\left[  \hat W_{m}(x-y)^p\right] \\
 &\le W_n(D)^p \sup_{y\in \bbZ^d}\bbE\left[  \hat W_{m}(y)^p\right].
\end{split}\end{equation}
and hence the desired result.
\end{proof}

To obtain \eqref{keystep} we  need to refine this approach a little bit. In the case $p\leq 2$, we split
\begin{align*}
	\bbE\left[ (\hat W_{n+m}(x))^p\right]=\bbE\left[ (\hat W_{n+m}(x))^2)^{p/2}\right]
\end{align*}
into two parts using a replica decomposition (the squared partition function is considered as an expectation for two independent random walks). The first part in this decomposition is the contribution to the partition function by pairs of trajectories that do not meet in the interval $(a_n,b_n)$, which yields the term  $\bbE[\hat W^{(a_n,b_n)}_{n}(x)^p ]$ after annealing (recall \eqref{annil}).
The remainder is the contribution to the partition function from pairs of trajectories that meet within the interval $(a_n,b_n)$, and to control it we rely on Lemma \ref{sympa}. More precisely, instead of using it for $D=\gO$ as in the proof of \eqref{rough2} we use it for a family of $D$ which has small probability under the polymer measure. This will give the small prefactor $\gep({\bf b} -{\bf a})$ in Proposition~\ref{keystep}. Here, we exploit the fact that when $({\bf b} -{\bf a})$ is small two replicas are unlikely to meet under the polymer measure as a consequence of the directed polymer CLT and invariance principle (expressed in \eqref{weakclt} and \eqref{invprincip}). Finally, the case $p>2$ is handled similarly except that we consider $\lceil p\rceil$ replicas instead of $2$.

\subsection{Proof of Theorem~\ref{factoth}}

We explain now how Theorem~\ref{factoth} is deduced from the results presented in the previous section. In the course of our proof, when $p<2$, we use a technical result which allows to compare the $L^p$ norm of a martingale increment in terms the increment of $L^p$ norm of the martingale. To explain the significance, note that when $p\ge 2$ and if
$X$ and $Y$ are two $L^p$-integrable random variables such that $\bbE[Y | X]=0$ and $X+Y\ge 0$ almost surely, then 
\begin{equation}\label{iiopi}
\bbE[|Y|^p]\le \bbE[ |Y+X|^p-|X|^p].
\end{equation}
Equality holds when $p=2$.
Indeed, by conditioning, we may restrict to the case where $X=a\ge 0$, and then we may check that the derivative of $\bbE[ (Y+a)^p- a^p ]$ is nonnegative using Jensen. The following provides a generalization for values $p\in (1,2)$.

\begin{lemma}{\cite[Lemma 7.1]{L25bb}}\label{otptilem}
 Let $p\in(1,2)$ and $X$ and $Y$ be two $L^p$-integrable random variables such that $\bbE[Y | X]=0$. Then, if $\bbE\left[ |Y+X|^p-|X|^p \right]\le \bbE[|X|^p]$, we have
 \begin{equation}\label{dfty}
 \bbE[|Y|^p]\le C_p \bbE[|X|^p]^{1-\frac{p}{2}}\bbE\left[ |Y+X|^p-|X|^p\right]^{p/2}
 \end{equation}
where $C_p>0$ is a universal constant that depends only on $p$.

\end{lemma}

\noindent Note that if $\bbE[|Y+X|^p-|X|^p]>\bbE[|X|^p]$, we obtain a similar conclusion as follows:
	\begin{align}\label{if_fails}
\bbE[|Y|^p] \leq 2^{p-1}\left(\bbE[|X+Y|^p]+\bbE[|X|^p]\right)  \leq 3\cdot 2^{p-1} \left(\bbE[|X+Y|^p]-\bbE[|X|^p]\right).
\end{align}
\noindent Now, our first step is to establish that for $p<p^*(\beta)$, the $p$-th moment of the point-to-point partition function is at most of order $n^{-\frac{d}{2}}$, that is
\begin{equation}\label{roughfirst}
 \sup_{x\in \bbZ^d}\bbE\left[(\hat W^{\beta}_{n}(x))^{p}\right]\le C_p n^{-\frac{dp}{2}}.
\end{equation}
The proof relies only on Proposition \ref{keystep}.
\begin{proof}[Proof of \eqref{roughfirst}]
Let us take  ${\bf a}=\frac{1}{2}-\delta$,  ${\bf b}=\frac{1}{2}$.
Using first \eqref{subtle} and then \eqref{superrough} for the second inequality we obtain
\begin{equation}\label{reccursion}\begin{split}
	\bbE\left[ \hat W_{2n}(x)^p\right] &\le \bbE\left[ \left(\hat W^{(a_{2n},b_{2n})}_{2n}\right)^p\right] +\gep(\delta)\sup_{y\in\bbZ^d}\bbE\left[ \hat W_{n}(y)^p\right]\\
					     &\le  C_p (2n\delta)^{-\frac{pd}{2}} +\gep(\delta)\sup_{y\in\bbZ^d}\bbE\left[ \hat W_{n}(y)^p\right].
 \end{split}
\end{equation}
We take $\delta>0$ small enough so that  $\gep(\delta)< 2^{-\frac{dp}{2}}$. Setting $n_k=2^k$
and
\begin{align*}
	w_k:= \sup_{x\in \bbZ^d} n_k^{dp/2} \bbE\left[\hat W_{n_k}(x)^p\right],
\end{align*}
we deduce from \eqref{reccursion} that
\begin{equation}
w_{k+1}\le 2^{dp/2}(C_{\delta,p} + \gep w_{k-1}),
\end{equation}
which implies that
$$\limsup_{k\to \infty} w_k \le  \frac{2^{dp/2}C_{\delta,p}}{1- 2^{dp/2}\gep}.$$
Finally for $n\in [n_k,n_{k+1})$ we have from \eqref{rough2}
 \begin{equation}
	 \max_{x\in \bbZ^d}\bbE[\hat W_{n}(x)^p]\le  \bbE[(W_\infty)^p] (n_k)^{-\frac{dp}{2}} a_k.\tag*{\qedhere}
\end{equation}
\end{proof}

\medskip Next, we  use \eqref{roughfirst} to prove \eqref{facto}.
First observe that it is sufficient to do it for one choice of $\ell_n$ and $m_n$. Indeed using independence if $\ell+m+1\le n$ and $\ell'\le \ell$ and $m'\le m$ we have
\begin{equation}\label{anylnmn}
 \| W^{\beta}_{\ell} \overleftarrow{W}^{\beta}_{m}(n,x)- W^{\beta}_{\ell'}\overleftarrow{W}^{\beta}_{m'}(n,x)    \|_p\\  \le
  \| W^{\beta}_{\ell}\|_p \|  W^{\beta}_{m}- W^{\beta}_{m'}\|_p +   \| W^{\beta}_{\ell}-W^{\beta}_{\ell'} \|_p \| W^{\beta}_{m'}\|_p
\end{equation}
Thus if two sequences $(\ell_n,m_n)$ and $(\ell'_n,m'_n)$ satisfy  the assumption of Theorem \ref{factoth}  then the convergence in $L^p$ of $W^{\beta}_n$ implies that for $p\in [1,p^*(\beta))$ we have
\begin{equation}\label{klcklc}
 \lim_{n\to \infty} \| W^{\beta}_{\ell_n} \overleftarrow{W}^{\beta}_{m_n}(n,x)- W^{\beta}_{\ell'_n} \overleftarrow{W}^{\beta}_{m'_n}(n,x)    \|_p=0.
\end{equation}
If $\ell'_n\le \ell_n$ and $m'_n\le m_n$, the convergence \eqref{klcklc} follows directly from \eqref{anylnmn}, if not we  insert the term $W^{\beta_c}_{\ell_n\wedge \ell'_n}\overleftarrow{W}^{\beta}_{m_n\wedge m'_n}(n,x)$ and use the triangle inequality before using \eqref{anylnmn}.
From the above, it is sufficient to prove the result with  $m_n=\ell_n=\lceil n/4\rceil$.

\medskip

\noindent Taking   ${\bf a}=\frac{1}{2}-\delta$,  ${\bf b}=\frac{1}{2}$ as before and combining \eqref{subtle} with \eqref{roughfirst} we obtain that
\begin{equation}\label{usethis}
\bbE[\hat W_{n}(x)^p]\le  \bbE\left[\hat W^{(a_n,b_n)}_{n}(x)^p \right]+  \gep(\delta)C_p (n/2)^{-\frac{dp}{2}}.
\end{equation}
When $p\ge 2$, \eqref{iiopi}  with
$X=\hat W^{(a_n,b_n)}_{n}(x)$ and $Y=\hat W_{n}(x)-\hat W^{(a_n,b_n)}_{n}(x)$   (recall \eqref{annil}) yields directly
\begin{equation}\label{eeeeee}
n^{\frac{dp}{2}} \bbE\left[\left|\hat W_{n}(x)-\hat W^{(a_n,b_n)}_{n}(x)\right|^p\right]\le   C\gep(\delta)=: \gep'(\delta).
\end{equation}
When $p<2$ we use instead Lemma \ref{otptilem}. Indeed, it the assumption of the lemma is satisfied, we obtain that for $n\ge n_0(\delta)$ sufficiently large
\begin{equation}\begin{split}\label{eeee}
n^{\frac{dp}{2}} \bbE\left[ \left|\hat W_{n}(x)-\hat W^{(a_n,b_n)}_{n}(x)\right|^p \right]& \le  \bbE\left[ n^{\frac{dp}{2}} \left|\hat W^{(a_n,b_n)}_{n}(x)\right|^p\right]^{1-\frac{p}{2}} \ \left(\gep(\delta)C_p 2^{\frac{dp}{2}}\right)^{p/2}\\
 &\le C'_p  \gep(\delta)^{p/2}=: \gep'(\delta),
 \end{split}
\end{equation}
where in the second line we made use of \eqref{lessrough}. Here it is important that $C_p'$ is independent of $\delta$, as observed after \eqref{lessrough}.
On the other hand, if the assumption of the lemma fails, \eqref{if_fails} and \eqref{usethis} also yield the final conclusion of \eqref{eeee} with a different constant $C_p'$.
Combining either \eqref{eeeeee} or \eqref{eeee} with Proposition \ref{factoslice} this yields
\begin{equation}
\limsup_{n\to \infty} n^{\frac{dp}{2}}\bbE\left[ \left|\hat W_{n}(x)-W^{\beta}_{a_n}\times \overleftarrow{W}^{\beta}_{n-b_n}(n,x) P(X_{n}=x)\right|^p \right]\le  \gep'(\delta).
\end{equation}
Finally using  to \eqref{anylnmn} yields
\begin{equation}
\limsup_{n\to \infty} n^{\frac{dp}{2}}\bbE\left[ \left|\hat W_{n}(x)-W^{\beta}_{\lceil n/4\rceil}\times \overleftarrow{W}^{\beta}_{\lceil n/4\rceil}(n,x) P(X_{n}=x)\right|^p \right]\le  \gep'(\delta).
\end{equation}
Since $\delta>0$ is arbitary and $\lim_{\delta \to 0}\gep'(\delta)=0$ this yields the desired conclusion.
\qed

\section{Proof of Proposition \ref{factoslice}}\label{sec_factoslice}

We need to prove that for
any $\bbZ^d$-valued sequence $(x_n)_{n\ge 0}$, we have
\begin{equation}
\limsup_{n\to \infty}n^{dp/2}\bbE\left[ \left|\hat W^{(a_n,b_n)}_{n}(x_n) -  W^{\beta}_{a_n}\times \overleftarrow{W}^{\beta}_{n-b_n}(n,x_n) P(X_{n}=x_n) \right|^p\right]=0.
\end{equation}
By extracting a subsequence,
we may assume without loss of generality that
 either $x_n/\sqrt{n}$ converges to ${\bf x}\in \bbR^d$ or $|x_n|/\sqrt{n}\to \infty$.
In the remainder of the proof we drop the superscript $n$ for $a_n,b_n$ and $x_n$ for better readability.

We let $\mu_n(z)=\frac{\hat W_n(z)}{W_n}$ denote the endpoint distribution of the directed polymer and for $m\le n$, we let $\overleftarrow \mu^{(n,x)}_m(z)$ denote the end point distribution associated to a polymer with partition function $\overleftarrow{W}^{\beta}_{m}(n,x)$, that is
$$\overleftarrow \mu^{(n,x)}_m(z):= \frac{E\left[ e^{\sum_{k=0}^m \beta\go_{n-k,x-X_k}-\gl(\beta)}\ind_{\{X_m=z\}}\right]}{\overleftarrow{W}^{\beta}_{m}(n,x)}. $$
We observe that the following identity holds (where ${\bf p}_k(z)=P(X_k=z)$)
\begin{equation}
\frac{
\hat W^{(a,b)}_{n}(x)}{ W^{\beta}_{a}\times \overleftarrow{W}^{\beta}_{n-b}(n,x)}
=\sum_{y,z\in \bbZ^d}\mu_a(y)\overleftarrow\mu^{(n,x)}_{n-b}(z){\bf p}_{b-a}(x+z-y).
\end{equation}
Now let us write the quantity whose expectation we must bound as a product
\begin{equation}
 \left|\hat W^{(a_n,b_n)}_{n}(x_n) -  W^{\beta}_{a_n}\times \overleftarrow{W}^{\beta}_{n-b_n}(n,x) P(X_{n}=x_n)\right|^p=: U_n \times (V_n)^p
\end{equation}
where
\begin{equation}\begin{split}
U_n&:=  \left(W^{\beta}_{a}\times \overleftarrow{W}^{\beta}_{n-b}(n,x)\right)^p,\\
V_n&:= n^{\frac{d}{2}}\left| \sum_{y,z\in \bbZ^d}\mu_{a}(y)\overleftarrow\mu^{(n,x)}_{n-b}(z){\bf p}_{b-a}(x+z-y) - {\bf p}_n(x)\right|.
\end{split}
\end{equation}
We need to prove that $\lim_{n\to \infty}
\bbE\left[ U_n \times (V_n)^p\right]=0$. We are going to show that
\begin{itemize}
 \item [(i)] The sequence $(U_n \times (V_n)^p)_{n \ge 1}$ is uniformly integrable.
 \item [(ii)] $\lim_{n\to \infty} V_n=0$ in probability.
\end{itemize}
We observe that $(U_n)_{n\ge 0}$ is uniformly integrable since for $q\in (p,p^*(\beta))$ we have by independence
\begin{equation}
  \bbE\left[ \left(W^{\beta}_{a}\times \overleftarrow{W}^{\beta}_{n-b}(n,x)\right)^q \right] =   \bbE\left[ \left(W^{\beta}_{a}\right)^q \right]\bbE\left[  \left(\overleftarrow{W}^{\beta}_{n-b}(n,x)\right)^q \right]
\end{equation}
and the r.h.s is uniformly bounded in $n$.
Furthermore, since $\liminf n^{-1}(b_n-a_n)>0$, we have
$$\esssup |V_n|\le \sup_{m\ge 1} m^{d/2}\left(\|{\bf p}_{b_m-a_m}\|_{\infty}+ \|{\bf p}_m\|_{\infty}\right)<\infty.$$
This concludes the proof of $(i)$. For $(ii)$ we use the local CLT (under $P$) to replace ${\bf p}$ by a Gaussian approximation,
\begin{equation}
 \frac{1}{2}(2\pi/d)^{d/2} V_n= \left|\left(\frac{n}{b-a}\right)^{d/2}\sum_{y,z\in \bbZ^d}\mu_{a}(y)\overleftarrow\mu^{(n,x)}_{n-b}(z) e^{-\frac{d |x+z-y|^2}{2(b-a)}} - e^{-\frac{d|x|^2}{2n}}\right|+o(1).
\end{equation}
In the case $x_n/n\to{\bf x}$, using also $a_n/n\to{\bf a}$ and $b_n/n\to{\bf b}$,  we further have
\begin{equation}
	\frac{1}{2}(2\pi/d)^{d/2} V_n= \left|\left(\frac{n}{b-a}\right)^{d/2}\sum_{y,z\in \bbZ^d}\mu_{a}(y)\overleftarrow\mu^{(n,x)}_{n-b}(z) e^{-\frac{d }{2({\bf b}-{\bf a})}\left|{\bf x} +\frac{z}{\sqrt n}-\frac{y}{\sqrt n}\right|^2} - e^{-\frac{d}{2}|{\bf x}|^2}\right|+o(1).
\end{equation}
Since $\mu_a(\cdot/\sqrt n)$ and $\overleftarrow \mu^{(n,x)}_{n-b}(\cdot/\sqrt n)$ converge, by the directed polymer CLT \eqref{weakclt}, to independent Gaussians with variances ${\bf a}$ and ${\bf 1-b}$, the first term in the above display converges to 
\begin{align*}
	\left(2\pi\sqrt{{\bf a}(1-{\bf b})({\bf b}-{\bf a})}\right)^{-d}\int_{\bbR^d\times \bbR^d}
 e^{-\frac{d|{\bf y}|^2}{2 {\bf a}}}  e^{-\frac{d|{\bf z}|^2}{2 (1-{\bf b})}}
e^{-\frac{d |{\bf x+z-y}|^2}{2({\bf b}-{\bf a})}}\dd{\bf y} \dd{\bf z}
=e^{-\frac{d|{\bf x}|^2}{2}}.
\end{align*}
If $\lim_{n\to \infty}|x_n|/\sqrt{n}=\infty$ the same CLT implies that
$$ \lim_{n\to \infty} \sum_{y,z\in \bbZ^d}\mu_{a}(y)\overleftarrow\mu^{(n,x)}_{n-b}(z) e^{-\frac{d |x+z-y|^2}{2(b-a)}}=0,$$
yielding the desired result in the two cases.
\qed

\section{Proof of Proposition \ref{keystep}}\label{sec_keystep}

We start with the case $p\le 2$ and then explain how to adapt the argument when $p> 2$.
\subsection{The case $p\le 2$} Like in the previous section, we frequently drop the superscript $n$ in $a_n$ and $b_n$ for better readability.
We let $\gO$ denote the set of random walk trajectories and we introduce the following subsets of $\gO^2$
\begin{equation}\begin{split}
 A= A(a,b,n,x):= \left\{ (X^{(1)},X^{(2)}) \ : \ \forall k\in (a,b), X^{(1)}_k\ne X^{(2)}_k ; X^{(1)}_{n}=X^{(2)}_{n}=x\right\},\\
  A'=A'(a,b,n,x):=\left\{(X^{(1)},X^{(2)}) \ : \ \exists k\in (a,b), X^{(1)}_k= X^{(2)}_k ; X^{(1)}_{n}=X^{(2)}_{n}=x\right\}.
 \end{split}
\end{equation}
We let  $W^{\otimes 2}_n(\cdot)$ denote the corresponding product measure on $\gO^2$ (dropping the dependence in $\beta$ in the notation for better readability).
Using $(a+b)^{\theta}\le a^{\theta}+b^{\theta}$ for $\theta\in (0,1)$, we have
\begin{equation}
 \bbE\left[\hat W^{\beta}_{n}(x)^p\right] \le \bbE\left[ W^{\otimes 2}_{n}(A)^{p/2} \right] + \bbE\left[ W^{\otimes 2}_{n}\left( A'\right)^{p/2} \right].
\end{equation}
Recalling \eqref{annil} and using the shorthand notation $\bbE_{(a,b)}:= \bbE\left[  \cdot \ | \ \cF_{\lint 1,a\rint \cup\lint b, n \rint}  \right]$, we observe that as a direct consequence of Fubini and the fact that there are no intersection in $(a,b)$,
\begin{equation}\begin{split}
               \bbE_{(a,b)}\left[W^{\otimes 2}_{n}(A)\right]
              & = E^{\otimes 2}\left[  \exp\left(\sum_{k\in \lint1,a\rint\cup\lint b,n\rint }\beta(\go_{k,X^{(1)}_k}+\go_{k,X^{(2)}_k}) -2\gl(\beta)\right)\ind_{\{ (X^{(1)},X^{(2)})\in A\}} \right]\\
            &\le \hat W^{(a,b)}_n(x)^2.
            \end{split}
\end{equation}
Hence using Jensen's inequality we have
\begin{equation}
               \bbE\left[ W^{\otimes 2}_{n}(A)^{p/2} \right] \le \bbE\left[
               \bbE_{(a,b)}\left[W^{\otimes 2}_{n}(A)\right]^{p/2} \right]\le
                \bbE\left[  \left( \hat W^{(a,b)}_n(x)\right)^p \right]
\end{equation}
which corresponds to the first term in \eqref{subtle}. Now it remains to show that  for large $n$
\begin{equation}\label{aprimm}
 \bbE\left[ W^{\otimes 2}_{n}\left( A'\right)^{p/2} \right]
\le \gep({\bf b}-{\bf a})\sup_{y\in \bbZ^d}\bbE\left[  \hat W_{n-b}(y)^p\right].
\end{equation}
We set $r= n^{1/2} ({\bf b}-{\bf a})^{1/4}$ and
define the events $B$ and $(C(z))_{z\in \bbZ^d}$ as follows

\begin{equation}\begin{split}\label{defBCz}
B_0&:=\left\{ X : \sup_{k\in[a,b]}\|X_k -X_a\|\ge r/2\right\}, \quad \quad  B:= B_0\cap\{X_n=x\}\\
C_0(z)&:=\{ X : X_a\in I(z)\} , \quad \quad \quad \quad \quad \quad \quad \  C(z):= C_0(z)\cap\{X_n=x\}
\end{split}\end{equation}
where in  $\|\cdot\|$ denote the $\ell^{\infty}$ norm on $\bbZ^d$ and $I(z):= rz+[-r,r)^d$ is a covering of $\bbR^d$ by overlapping cubes of sidelength $2r$.
We further set
\begin{equation}
 \overline B^{(i)}:=\{X^{(i)}_n\in B_0\} \cap \{  X^{(1)}_n=X^{(2)}_n=x \} \quad
 \text{ and }  \quad \overline C(z)=:C(z)\times C(z),
\end{equation}
and observe that
\begin{equation}\label{incluz}
 A'\subset  \overline B^{(1)} \cup \overline B^{(2)}\cup \bigcup_{z\in \bbZ^d}  \overline C(z).
\end{equation}
Indeed if $X^{(1)}$ and $X^{(2)}$ do not belong to a common $C(z)$, then $\|X^{(1)}_a-X^{(2)}_a\|> r$ and at least one of them needs to travel a distance $r/2$ before $b$ to create an intersection.
The inclusion \eqref{incluz} expresses the fact that if two trajectories belong to $A'$ either they are close to each other at time $a_n$ or they have a slightly superdiffusive behavior in the time interval $(a_n,b_n)$, at least when $({\bf b}-{\bf a})$ is small.
Using subadditivity and simplifying the expression (using the fact that the events are products) yields
\begin{equation}\begin{split}\label{trapiz}
 \bbE\left[ W^{\otimes 2}_{n}\left( A'\right)^{p/2} \right]&\le   \bbE\left[ W^{\otimes 2}_n( \overline B^{(1)})^{p/2}\right]+\bbE\left[ W^{\otimes 2}_n( \overline B^{(2)})^{p/2}\right]+ \sum_{z\in \bbZ^d}\bbE\left[ W^{\otimes 2}_n(\overline C(z))^{p/2}\right]\\
 &\le 2\bbE\left[  W_n(B)^{p/2} \hat W_n(x)^{p/2}\right]+  \sum_{z\in \bbZ^d}\bbE\left[ W_n(C(z))^{p}\right]\\
 &\le 2\bbE\left[  W_n(B)^p\right]^{1/2} \bbE\left[ \hat W_n(x)^p\right]^{1/2}+  \sum_{z\in \bbZ^d}\bbE\left[ W_n(C(z))^{p}\right]
 \end{split} \end{equation}
Recall \eqref{Wasmeasure} and note that both terms in the r.h.s.\ of \eqref{trapiz} are of the form $\bbE\left[ \hat W_n(D,x)^{p}\right]$ for some event $D$ defined in terms of the first $b$ steps of $X$.
Applying this to $B_0$ and $C_0(z)$ and $\gO$ in \eqref{trapiz},
Lemma~\ref{sympa} yields
\begin{equation}
 \bbE\left[ W^{\otimes 2}_{n}\left( A'\right)^{p/2} \right]
 \le  \left( 2\bbE\left[ W_b(B_0)^p\right]^{1/2}\bbE\left[ (W_b)^p\right]^{1/2} + \sum_{z\in \bbZ^d}\bbE\left[W_b(C_0(z))^p \right] \right)\sup_{y\in \bbZ^d}\bbE\left[  \hat W_{n-b}(y)^p\right]
\end{equation}
and thus to conclude the proof of \eqref{aprimm} we are going to show that for $n$ sufficiently large
\begin{equation}\label{tbbprov}
 2\bbE\left[ W_b(B_0)^p\right]^{1/2}\bbE\left[ (W_b)^p\right]^{1/2} + \sum_{z\in \bbZ^d}\bbE\left[W_b(C_0(z))^p \right]  \le C \left({\bf b}-{\bf a}\right)^{\frac{d(p-1)}{4}}.
\end{equation}
Using the invariance principle \eqref{invprincip} for the directed polymer (combined with a Portmanteau type argument since here $\varphi$ is the indicator of the set whose boundary is negligible under the Wiener measure) we have the following convergence in probability for ${\bf r}=({\bf b}-{\bf a})^{1/4}$
\begin{equation}
 \lim_{n\to \infty} W^{\beta}_n(B_0)= W_\infty \bP\left( \sup_{t\in[{\bf a},{\bf b}]} |B_t-B_{\bf a}|\ge  {\bf r}/2 \right).
\end{equation}
Hence, using dominated convergence under $\bbP$ (by Doob's inequality $\sup_{n\ge 1} W_n^p\in L^1$),
we have
\begin{equation}\label{leBB}
 \lim_{n\to \infty} \bbE\left[ W^{\beta}_b(B_0)^p\right]=
 \bbE\left[(W^{\beta}_{\infty})^p\right] \bP\left( \sup_{t\in[{\bf a},{\bf b}]} |B_t-B_{\bf a}|\ge  {\bf r}/2 \right)^p.
\end{equation}
In the same manner we obtain that
\begin{equation}\label{leCC}
  \lim_{n\to \infty} \bbE\left[ W^{\beta}_b(C_0(z))^p\right]
  = \bbE\left[(W^{\beta}_{\infty})^p\right] \bP\left( B_{{\bf a} }\in [-{\bf r},{\bf r})^d+{\bf r} z \right)^p.
\end{equation}
To to take the sum over $z\in \bbZ^d$, we must make sure that
\begin{equation}\label{TOUT}
 \lim_{R\to \infty} \limsup_{n\to \infty}\sum_{\|z\|>R}\bbE\left[ W^{\beta}_b(C_0(z))^p\right]=0.
\end{equation}
Using $\sum a_i^{p}\le (\sum a_i)^p$  for $a_i\geq 0$ and $p\geq 1$ and $\sum_{z\in \bbZ^d} \ind_{C_0(z)}\le 2^d$, we obtain
\begin{equation}
  \sum_{\|z\|>R} \bbE\left[ W^{\beta}_b(C_0(z))^p\right]
  \le  2^{dp} \bbE\left[ W^{\beta}_b\left( \bigcup_{\|z\|>R}C_0(z)\right)^p\right]
\end{equation}
and we may conclude that \eqref{TOUT} holds. Indeed, as in the case of \eqref{leCC} and \eqref{leBB}, the invariance principle for directed polymers implies that
\begin{equation}
 \lim_{n\to \infty}\bbE\left[ W^{\beta}_b\left( \bigcup_{\|z\|>R}C_0(z)\right)^p\right]= \bbE\left[(W^{\beta}_{\infty})^p\right] \bP\big(\| B_{{\bf a} } \| \ge (R-1){\bf r} \big)^p.
\end{equation}
and hence that \eqref{TOUT} holds. Combining \eqref{leBB}, \eqref{leCC} and \eqref{TOUT} we obtain
\begin{multline}\label{fatlimit}
\lim_{n\to \infty} \left(2\bbE\left[ W_b(B_0)^p\right]^{1/2}\bbE\left[ (W_b)^p\right]^{1/2} + \sum_{z\in \bbZ^d}\bbE\left[W_b(C_0(z))^p \right]\right) \\
=\bbE\left[ (W_\infty)^p\right] \left(2\bP\left( \sup_{t\in[{\bf a},{\bf b}]} |B_t-B_{\bf a}|\ge  {\bf r} \right)^{\frac p 2}+  \sum_{z\in \bbZ^d} \bP\left( B_{\bf a}\in [-{\bf r},{\bf r}) -z{\bf r}\right)^p\right)
\le C {\bf r}^{\frac{d}{2}(p-1)}.
 \end{multline}
The last inequality uses the fact that ${\bf a}\ge 1/4$ to estimate the second term (since ${\bf r}= ({\bf b}- {\bf a})^{1/4}$ the first term is of order $\exp(-c {\bf r}^{-2})$). This concludes the proof of \eqref{tbbprov} and hence of the Proposition.\qed

\subsection{The case $p>2$}\label{p_larger_2}
 Let us assume that $p\in (q-1,q]$ for an integer $q\ge 3$.
 We introduce  $q$-dimensionals analogue of $A$ and $A'$,
 \begin{equation}\begin{split}
  A:= \left\{ \forall k\in (a,b)\ \forall i\ne j, \  X^{(i)}_k\ne X^{(j)}_k\} \cap   \{ \forall i \in \lint 1,q\rint,  X^{(i)}_{n}=x\right\},\\
    A':= \left\{ \exists k\in (a,b)\ \exists i\ne j, \  X^{(i)}_k= X^{(j)}_k\} \cap   \{ \forall i \in \lint 1,q\rint,  X^{(i)}_{n}=x\right\},
    \end{split}
\end{equation}
and observe that
\begin{equation}
 \bbE\left[\hat W^{\beta}_{n}(x)^p\right] \le \bbE\left[ W^{\otimes q}_{n}(A)^{p/q} \right] + \bbE\left[ W^{\otimes q}_{n}\left( A'\right)^{p/q} \right].
\end{equation}
The same computation as in the $q=2$ case yields
\begin{equation}
  \bbE\left[ W^{\otimes q}_{n}(A)^{p/q} \right]\le \bbE\left[\hat W^{(a,b)}_
  {n}(x)^p\right]
\end{equation}
To bound   $\bbE\left[ W^{\otimes q}_{n}(A')^{p/q} \right]$
we recall  the definition of $B_0,B, C_0(z)$ and $C(z)$ in \eqref{defBCz}
and set
\begin{equation}\begin{split}
 \overline B^{(i)}&:= \{ X^{(i)}\in B_0 \} \cap\{\forall k\in \lint 1, q\rint,\ X^{(k)}_n=x\},\\
\overline C^{(i,j)}(z)&:= \{ X^{(i,j)}\in C_0(z) \}\cap\{\forall k\in \lint 1, q\rint,\ X^{(k)}_n=x\}.
 \end{split}\end{equation}
We observe that (for the same reason as in the $q= 2$ case)
$$A' \subset \bigcup_{i=1}^q\overline B^{(i)}  \cup \uniontwo{z\in \bbZ^d}{1\le i<j\le q} \overline C^{(i,j)}_n(z).$$
Using subbaditivity symetries and the fact that the events are cylinders, we deduce from the above that
\begin{multline}
 \bbE\left[ W^{\otimes q}_{n}\left( A'\right)^{p/q} \right]
 \le  q \bbE\left[ \hat W^{\beta}_n(x)^{\frac{(q-1)p}{q}}  W_n (B)^{\frac{p}{q}} \right] + \frac{q(q-1)}{2}\sum_{z\in \bbZ^d}
 \bbE\left[ \hat W^{\beta}_n(x)^{\frac{(q-2)p}{q}}  W_n (C(z))^{\frac{2p}{q}} \right]\\
 \le q \bbE\left[ \hat W^{\beta}_n(x)^{p}\right]^{\frac{q-1}{2}} \bbE\left[  W_n (B)^{p} \right]^{\frac 1 q} + \frac{q(q-1)}{2}\sum_{z\in \bbZ^d}
 \bbE\left[ \hat W^{\beta}_n(x)^p\right]^{\frac{(q-2)}{q}}\bbE\left[ W_n (C(z))^{\frac{2p}{q}} \right]^{\frac 2 q}.
\end{multline}
Using Lemma \ref{sympa}, we have
\begin{multline}
  \bbE\left[ W^{\otimes q}_{n}\left( A'\right)^{p/q} \right]
\leq \sup_{y\in \bbZ^d} \bbE\left[ \hat W^{\beta}_{n-b}(y)^p\right]\\
\times
\left( q\bbE\left[  (W^{\beta}_b)^{p}\right]^{\frac{q-1}{2}} \bbE\left[  W_b (B_0)^{p} \right]^{1/q}+ \frac{q(q-1)}{2}\sum_{z\in \bbZ^d}
 \bbE\left[ (W^{\beta}_b)^p\right]^{\frac{(q-2)}{q}}\bbE\left[ W_b (C_0(z))^{p} \right]^{2/q}\right).
\end{multline}
Then we may conclude by observing that for the same reason as \eqref{fatlimit}, the second line in the above equation converges to
\begin{equation}
 \bbE\left[ (W_\infty)^p\right] \left(q\bP\left( \sup_{t\in[{\bf a},{\bf b}]} |B_t-B_{\bf a}|\ge  {\bf r} \right)^{\frac p q}\!\! + \frac{q(q-1)}{2} \sum_{z\in \bbZ^d} \bP\left( B_{\bf a}\in [-{\bf r},{\bf r}) +z{\bf r}\right)^\frac{2p}{q}\right)
\le C {\bf r}^{d(\frac{2p}{q}-1)}.
\end{equation}
To ensure that the limit can be exchanged with the sum, we replace \eqref{TOUT} by
\begin{align*}
	\lim_{R\to\infty}\limsup_{n\to\infty} \sum_{\|z\|>R}\E[W_b^\beta(C_0(z))^p]^{2/q}\leq\lim_{R\to\infty}\limsup_{n\to\infty} \sum_{\|z\|>R}P(C_0(z))^{\frac{2(p'-p)}{q(p'-1)}}=0,
\end{align*}
for some $p'\in(p,\p)$.

\qed

\section{Proof of Theorem \ref{condinvpr}}\label{sec_condinvpr}
We let  $\hat P^{\beta,\go}_{n,x}:= P^{\beta,\go}_{n}(  \cdot \ | \ X_n=x)$ denote  the conditional measure
\begin{equation}
 \hat E^{\beta,\go}_{n,x}[g(X)]:= \frac{1}{\hat W^{\beta}_{n}(x)}  E\left[ \exp\left( \sum_{k=1}^n \beta\go_{k,X_k}-\gl(\beta)\right)\ind_{\{X_n=x\}}g(X)\right]
\end{equation}
and simply write
$\hat P_{n,x}$ when $\beta=0$.
Since we already know that $\lim_{n\to \infty}\hat E_{n,x_n}[\varphi(X^{(n)})]|=\bE\left[\varphi(B^{({\bf x})})\right]$ we need to show that for every continous $\varphi\colon\cC\to [0,1]$ we have
\begin{equation}
 \lim_{n\to \infty}\left|\hat E^{\beta,\go}_{n,x_n}\left[\varphi(X^{(n)})\right]- \hat E_{n,x_n}\left[\varphi(X^{(n)})\right]\right|=0.
\end{equation}
We set $m_n=n^{1/8}$ and we let $\hat P^{\beta,(1)}_{n,x}$ denote the polymer measure with an environment only present within distance $m_n$ of the two endpoints of the polymer.
\begin{equation}
\hat E^{\beta,(1)}_{n,x}[g(X)] := \frac{1}{\hat W^{\beta,(1)}_{n}(x)}  E\left[ \exp\left( \sum_{k\in \lint 1, m_n\rint \cup \lint n-m_n,n\rint} \beta\go_{k,X_k}-\gl(\beta)\right)\ind_{\{X_n=x\}}g(X)\right].
\end{equation}
where $\hat W^{\beta,(1)}_{n}(x)= \hat W^{\beta,(m_n,n-m_n)}_{n}(x)$ with the terminology introduced in \eqref{defwaabb}.
 We have
\begin{multline}\label{moultix}
 \left|\hat E^{\beta,\go}_{n,x}[\varphi(X^{(n)})]- \hat E_{n,x}[\varphi(X^{(n)})]\right|\\
 \le
 \left|\hat E^{\beta,\go}_{n,x}[\varphi(X^{(n)})]- \hat E^{\beta,(1)}_{n,x}[\varphi(X^{(n)})]\right|+ \left|\hat E^{\beta,(1)}_{n,x}[\varphi(X^{(n)})]- \hat E_{n,x}[\varphi(X^{(n)})] \right|
\end{multline}
We need to show that both terms on the r.h.s.\ vanish to zero.
For the second term we observe that under $\hat P^{\beta,(1)}_{n,x}$,
the restriction of $X$ to the interval $[m_n,n-m_n]$ is a random walk bridge whose starting and ending point are random but at most at a distance $m_n$ from $0$ and $x_n$ respectively. Since $m_n=o(\sqrt{n})$ the first $m_n$ and last $m_n$ step do not affect the scaling limit. For completeness we perform this routine verification in the form of the following lemma, which allows to couple $\hat P^{\beta,(1)}_{n,x}$ with $\hat P_{n,x}$ on a smaller interval $[\ell_n,n-\ell_n]$ (our choice of $\ell_n$ and $m_n$ are made for convenience to make the verification process easier).

\begin{lemma}\label{neuf}
Let $\ell_n:=n^{1/3}$. Given $x_n$ such that $|x_n|/\sqrt{n}$ is bounded we have
\begin{equation}
 \lim_{n\to \infty}  \left\|\hat P^{\beta,(1)}_{n,x_n}\left( (X_k)_{k\in \lint \ell_n, n-\ell_n\rint}\in \cdot \right)- P_{n,x}\left( (X_k)_{k\in \lint \ell_n, n-\ell_n\rint}\in \cdot \right)\right\|_{TV}=0.
\end{equation}
\end{lemma}

The above result entails the existence of a coupling of trajectories $X$ with distribution $\hat P^{\beta,(1)}_{n,x}$ and  $Y$ with distribution $\hat P_{n,x}$ such that
with probability tending to $1$, $|X^{(n)}_t-Y^{(n)}_t|\le \ell_n/\sqrt{n}$ for all $t\in\lint 1,n\rint$. Indeed, the trajectories coincide in the interval $\lint \ell_n, n-\ell_n\rint$ with probability $1-o(1)$ and can only move distance $\ell_n$ in the remaining time. Thus  we have
\begin{equation}
 \lim_{n\to \infty}\left|\hat E^{\beta,(1)}_{n,x}[\varphi(X^{(n)})]- \hat E_{n,x}[\varphi(X^{(n)})] \right|=0.
\end{equation}

\begin{proof}[Proof of Lemma \ref{neuf}]
Using the spatial Markov property, it is sufficient to show that the joint distributions of $(X_{\ell},X_{n-\ell})$ are close under
 $\hat P^{\beta,(1)}_{n,x_n}$ and $\hat P_{n,x_n}$.
We have
\begin{multline}
 \hat P^{\beta,(1)}_{n,x}( (X_{\ell},X_{n-\ell})=(y,z)) \\ = \frac{1}{\cZ}\sum_{y_1,z_1\in \bbZ^d} \mu(y_1) P(X_{\ell-m}=y-y_1)P(X_{n-2\ell}=z-y)P(X_{\ell-m}=z_1-z)   \nu(z_1)
\end{multline}
where $\cZ$ is an appropriate renormalization factor and $\mu$ and  $\nu$ correspond to the end-point distribution of some polymer measure (the important point is that their support is included in a ball of radius $m_n$ around the origin and $x_n$ respectively).
Since $P(X_{n-2\ell}=z-y)$ does not vary much with the choice of $z$ and $y$, provided that $|y|\le \ell$ and $|z-x|\le \ell$, we have
\begin{equation}
\lim_{n\to \infty} \sup_{z,y \  : \   \hat P_{n,x}( (X_{\ell},X_{n-\ell})=(y,z))>0}   \ \left|\frac{P(X_{n-2\ell}=z-y)}{\cZ}-1\right|=1
\end{equation}
which implies
\begin{equation}\begin{split}
 \hat P^{\beta,(1)}_{n,x}&( (X_{\ell},X_{n-\ell})=(y,z))\\
 &=
 \sum_{y_1\in \bbZ^d} \mu(y_1) P(X_{\ell-m}=y-y_1)
 \sum_{z_2 \in \bbZ^d} \nu(z_1) P(X_{\ell-m}=z_1-z)(1+\gep(y,z,n))
 \end{split}
\end{equation}
with $\lim_{n\to \infty} \sup_{z,y}|\gep(y,z,n)|= 0$. Note that similarly we have
\begin{equation}
 \hat P_{n,x}( (X_{\ell},X_{n-\ell})=(y,z))
 =P(X_{\ell}=y)
 P(X_{\ell}=x-z)(1+\gep'(y,z,n)).
\end{equation}
We have approximated both  $\hat P^{\beta,(1)}_{n,x}$ and $\hat P_{n,x}( (X_{\ell},X_{n-\ell})\in \cdot )$ by product measures.
To conclude that the two distribution are close in total variation we need to check that
\begin{equation}\begin{split}
\lim_{n\to \infty}\sum_{y\in \bbZ^d} \bigg|\sum_{y_1\in \bbZ^d} \mu(y_1) P(X_{\ell-m}=y-y_1)- P(X_{\ell}=y)\bigg|=0,\\
\lim_{n\to \infty}\sum_{z\in \bbZ^d} \bigg|\sum_{z_1\in \bbZ^d} \nu(y_1) P(X_{\ell-m}=z_1-z)- P(X_{\ell}=x_n-z)\bigg|=0,
\end{split}
\end{equation}
which since thanks to our scale separation, is a consequence of the local CLT for the simple random walk.

\end{proof}

To deal with the second term in \eqref{moultix}, we recall the notation \eqref{Wasmeasure} (which we extend to $\hat W^{\beta,(1)}$) and observe that
\begin{multline}\label{preskfini}
  \left|\hat E^{\beta,\go}_{n,x}\left[\varphi(X^{(n)})\right]- \hat E^{\beta,(1)}_{n,x}\left[\varphi(X^{(n)})\right]\right|=\left| \frac{\hat W^{\beta}_{n}(\varphi(X^{(n)}),x)}{\hat W^{\beta}_{n}(x)}- \frac{\hat W^{\beta,(1)}_{n}(\varphi(X^{(n)}),x)}{\hat W^{\beta,(1)}_{n}(x)} \right|\\
  \le \left| \frac{\hat W^{\beta}_{n}(\varphi(X^{(n)}),x)-\hat W^{\beta,(1)}_{n}(\varphi(X^{(n)}),x)}{\hat W^{\beta}_{n}(x)} \right| +
  \left| \frac{\hat W^{\beta,(1)}_{n}(\varphi(X^{(n)}),x)}{\hat W^{\beta}_{n}(x)}- \frac{\hat W^{\beta,(1)}_{n}(\varphi(X^{(n)}),x)}{\hat W^{\beta,(1)}_{n}(x)} \right|
  \\
  \le \frac{1}{\hat W^{\beta}_{n}(x)}\left(\left|\hat W^{\beta}_{n}(\varphi(X^{(n)}),x)- \hat W^{\beta,(1)}_{n}(\varphi(X^{(n)}),x)\right| + \left|\hat W^{\beta,(1)}_{n}(x)-\hat W^{\beta}_{n}(x)\right|\right)
\end{multline}
We conclude by showing that for every $\varphi\colon\cC\to [0,1]$ we have
\begin{equation}\begin{split}\label{twothings}
  \lim_{n\to \infty}n^{d/2}\bbE\left[\left|\hat W^{\beta,(1)}_{n}(\varphi(X^{(n)},x_n)- \hat W^{\beta}_{n}(\varphi(X^{(n)},x_n)\right|\right]&=0,\\
  \lim_{\gep \to 0}\limsup_{n\to \infty} \bbP\left( n^{d/2} \hat W^{\beta}_{n}(x_n)>\gep \right)&=1.
  \end{split}
\end{equation}
The second line implies that the prefactor in \eqref{preskfini} is tight while the first line (applied to $\varphi\equiv 1$ to control the second term) implies that the content of the bracket converges to zero.

\medskip

\noindent To establish the first line of \eqref{twothings} we need to prove  the following for any $g\colon\gO\to [0,1]$,
\begin{equation}
 \bbE\left[\left|\hat W^{\beta,(1)}_{n}(g,x_n)- \hat W^{\beta}_{n}(g,x_n)\right| \right]
\le \bbE\left[ \left|\hat W^{\beta,(1)}_{n}(x_n)- \hat W^{\beta}_{n}(x_n)\right|\right].
 \end{equation}
This is done by replicating the proof of \cite[Lemma 2.1]{L25}. Moreover, as a consequence of \eqref{facto} (for $p=1$) we have
\begin{equation}
\lim_{n\to \infty}n^{d/2}\bbE\left[\left |\hat W^{\beta,(1)}_{n}(x_n)- \hat W^{\beta}_{n}(x_n)\right|\right]=0,
\end{equation}
since $\hat W^{\beta,(1)}_{n}(x_n)$ is very well approximated by the product $W^{\beta}_{m_n}\times \overleftarrow{W}^{\beta}_{m_n}(n,x_n)\times P(X_{n}=x_n)$.
For the second line in \eqref{twothings}, we observe that from \eqref{facto} it is sufficient to prove that
\begin{equation}
   \lim_{\gep \to 0}\limsup_{n\to \infty} \bbP\left[ W^{\beta}_{m_n}\times \overleftarrow{W}^{\beta}_{m_n}(n,x_n)>\gep \right]=1.
\end{equation}
which follows immediately from the fact that both
$W^{\beta}_{m_n}$ and $\overleftarrow{W}^{\beta}_{m_n}(n,x_n)$ converge in distribution to a positive random variable. \qed

\subsection*{Acknowledgements} S.J. was supported by JSPS Grant-in-Aid for Scientific Research 23K12984 and 25K00911. {H.L. was supported by a productivity grant from CNPq, a CNE grant from FAPERj and by funding from the European Union’s Horizon Europe
research and innovation programme under the Marie Skłodowska-Curie grant agreement No 101197621.
}

\bibliographystyle{plain}
\bibliography{ref.bib}

\begin{thebibliography}{10}

\bibitem{AZ96}
Sergio Albeverio and Xian~Yin Zhou.
\newblock A martingale approach to directed polymers in a random environment.
\newblock {\em J. Theoret. Probab.}, 9(1):171--189, 1996.

\bibitem{BT10}
Quentin Berger and Fabio~L. Toninelli.
\newblock On the critical point of the random walk pinning model in dimension
  d=3.
\newblock {\em Electronic J. Probab.}, 15:654--683, 2010.

\bibitem{BS10}
Matthias Birkner and Rongfeng Sun.
\newblock Annealed vs quenched critical points for a random walk pinning model.
\newblock {\em Ann. Inst. Henri Poincar{\'e}, Probab. Stat.}, 46(2):414--441,
  2010.

\bibitem{BS11}
Matthias Birkner and Rongfeng Sun.
\newblock Disorder relevance for the random walk pinning model in dimension 3.
\newblock {\em Ann. Inst. Henri Poincar{\'e}, Probab. Stat.}, 47(1):259--293,
  2011.

\bibitem{B89}
Erwin Bolthausen.
\newblock A note on the diffusion of directed polymers in a random environment.
\newblock {\em Comm. Math. Phys.}, 123(4):529--534, 1989.

\bibitem{CH02}
Philippe Carmona and Yueyun Hu.
\newblock On the partition function of a directed polymer in a {G}aussian
  random environment.
\newblock {\em Probab. Theory Related Fields}, 124(3):431--457, 2002.

\bibitem{CSY03}
Francis Comets, Tokuzo Shiga, and Nobuo Yoshida.
\newblock Directed polymers in a random environment: path localization and
  strong disorder.
\newblock {\em Bernoulli}, 9(4):705--723, 2003.

\bibitem{CY06}
Francis Comets and Nobuo Yoshida.
\newblock Directed polymers in random environment are diffusive at weak
  disorder.
\newblock {\em Ann. Probab.}, 34(5):1746--1770, 2006.

\bibitem{HKL23}
Tobias Hurth, Konstantin Khanin, Beatriz Navarro~Lameda, and Fedor Nazarov.
\newblock On a factorization formula for the partition function of directed
  polymers.
\newblock {\em J. Stat. Phys.}, 190(10):Paper No. 165, 32, 2023.

\bibitem{IS88}
John~Z. Imbrie and Thomas Spencer.
\newblock Diffusion of directed polymers in a random environment.
\newblock {\em J. Statist. Phys.}, 52(3-4):609--626, 1988.

\bibitem{J22}
Stefan Junk.
\newblock Fluctuations of partition functions of directed polymers in weak
  disorder beyond the {L}2-phase.
\newblock {\em Ann. Probab.}, 53(2):--, 2025.

\bibitem{J_JEMS}
Stefan Junk.
\newblock Local limit theorem for directed polymers beyond the $l^{2}$-phase.
\newblock {\em Journal of the European Mathematical Society}, May 2025.

\bibitem{JL_in_prep}
Stefan Junk and Hubert Lacoin.
\newblock Near-critical behavior of directed polymers in dimensions $d\ge 3$.
\newblock in preparation.

\bibitem{JL_ASENS}
Stefan {Junk} and Hubert {Lacoin}.
\newblock {Strong disorder and very strong disorder are equivalent for directed
  polymers}.
\newblock {\em arXiv e-prints}, page arXiv:2402.02562, February 2024.
\newblock to appear in \emph{Ann. Sci. Éc. Norm. Supér.}

\bibitem{JL_ORB}
Stefan Junk and Hubert Lacoin.
\newblock Coincidence of critical points for directed polymers for general
  environments and random walks.
\newblock {\em Orbita Math.}, 3(1):89--125, 2026.

\bibitem{L25bb}
Hubert {Lacoin}.
\newblock {The localization transition for the directed polymer in a random
  environment is smooth}.
\newblock {\em to appear in Pub. Math. IHES}.

\bibitem{L25}
Hubert Lacoin.
\newblock A short proof of diffusivity for the directed polymers in the weak
  disorder phase.
\newblock {\em Electron. Commun. Probab.}, 30:--, 2025.

\bibitem{LL10}
Gregory~F. Lawler and Vlada Limic.
\newblock {\em Random walk: a modern introduction}, volume 123 of {\em
  Cambridge Studies in Advanced Mathematics}.
\newblock Cambridge University Press, Cambridge, 2010.

\bibitem{sinai95}
Yakov~G. Sinai.
\newblock A remark concerning random walks with random potentials.
\newblock {\em Fund. Math.}, 147(2):173--180, 1995.

\bibitem{V06}
Vincent Vargas.
\newblock A local limit theorem for directed polymers in random media: the
  continuous and the discrete case.
\newblock {\em Ann. Inst. H. Poincar\'e{} Probab. Statist.}, 42(5):521--534,
  2006.

\end{thebibliography}

\end{document}